\documentclass[a4paper,11pt]{article}
\usepackage[T2A]{fontenc}
\usepackage[utf8]{inputenc}
\usepackage[russian,english]{babel}

\usepackage{latexsym}
\usepackage{geometry}
\usepackage{graphicx}
\usepackage{amsfonts}
\usepackage{amsthm}
\usepackage{amsmath}
\usepackage{amssymb}
\usepackage{xcolor}

\usepackage{accents}
\usepackage{indentfirst}
\usepackage{float}
\usepackage{mathtools}
\usepackage[noadjust]{cite}
\usepackage{array}
\usepackage{hyperref}

\usepackage{enumitem}
\setlist[enumerate]{font=\upshape,label=(\arabic*)}
\setlist{itemsep=1pt,topsep=4pt}

\newtheorem{theorem}{Theorem}[section]
\newtheorem{corollary}[theorem]{Corollary}
\newtheorem{proposition}[theorem]{Proposition}

\theoremstyle{definition}
\newtheorem{definition}[theorem]{Definition}
\newtheorem{notation}[theorem]{Notation}
\newtheorem{example}[theorem]{Example}
\newtheorem{remark}[theorem]{Remark}

\numberwithin{equation}{section}
\allowdisplaybreaks

\def\N{\mathbb{N}}
\def\Z{\mathbb{Z}}

\def\CD{{\rm CD}}
\def\n{{\rm n}}
\def\t{{\rm t}}

\providecommand{\keywords}[1]{\textbf{Keywords:} #1}
\providecommand{\msc}[1]{\textbf{MSC 2020:} #1}

\newcommand\freefootnote[1]{%
    \bgroup
    \renewcommand\thefootnote{\fnsymbol{footnote}}%
    \renewcommand\thempfootnote{\fnsymbol{mpfootnote}}%
    \footnotetext[0]{#1}%
    \egroup
}

\begin{document}

\title{Embedding of Barnes--Wall lattices\\
into Cayley--Dickson algebras}
\author{
Svetlana Zhilina$^{a}$
}
\date{\small \em
$^a$Lomonosov Moscow State University, Moscow, 119991, Russia
}

\freefootnote{This work was supported by the Russian Science Foundation (project No.~26-11-00001).}

\freefootnote{Email address: \texttt{s.a.zhilina@gmail.com}}

\maketitle

\begin{abstract}
The family of Barnes--Wall lattices includes the balanced Barnes--Wall lattice~$M_n$ over $\Z[\sqrt{2}]$ and two ordinary Barnes--Wall lattices~$L_n$ and~$L'_n$ over~$\Z$, which are geometrically equivalent to each other. Embedding of these lattices into the Cayley--Dickson algebras of an appropriate dimension is considered. It is shown that for all $n \geq 2$ the lattices~$M_n$ and~$L_n$ are closed under multiplication and are generated by the sets of pairwise products of their elements. The lattices $L_n \cdot L'_n$, $L'_n \cdot L_n$, and $L'_n \cdot L'_n$, generated by pairwise products of elements of the lattices~$L_n$ and~$L'_n$, also are studied.
\end{abstract}


\msc{11H06, 17A20, 17A45}

\section{Introduction}

The problem of dense packing of equal spheres goes back to Lagrange in the two-dimensional case and to Gauss in the three-dimensional one. According to the result of Viazovska~\cite{Viazovska8}, the $E_8$ lattice realizes the densest sphere packing in the eight-dimensional space. Later Viazovska and her coauthors~\cite{Viazovska24} proved a similar statement for the Leech lattice $\Lambda_{24}$ of dimension~$24$. In other dimensions greater than three this problem still remains open, despite the active work of many leading scientists in this direction.

It is proved that, in the four-dimensional space, the $D_4$ lattice gives the densest packing of equal spheres among those packings such that the centers of the spheres form a lattice~\cite[Chapter~1, p.~10]{ConwaySloane}. However, according to the currently known data, in an arbitrary dimension the densest packing need not be a lattice packing: for example, in dimensions~$10$, $11$, and~$13$ the packings~$P_{10c}$, $P_{11a}$, and~$P_{13a}$ turn out to be more efficient than all {\em known} lattice packings~\cite[Chapter~5, p.~140]{ConwaySloane}. In dimension~$16$ the best known result is given by the $16$-dimensional Barnes--Wall lattice. It can be seen from~\cite[Chapter~1, p.~14, Figure~1.5]{ConwaySloane} that it attains a local peak of the sphere packing density which brings it closer to the Rogers' upper bound.

In general, the Barnes--Wall lattices form a family of $2^n$-dimensional lattices, $n \geq 0$, which was first introduced in~\cite{BarnesWall}. For each~$n$ this family contains two geometrically equivalent lattices nested in each other, and the first members of this sequence include~$\Z$, $\Z^2$, $D_4$, and~$E_8$ for $n = 0$, $1$, $2$, and~$3$, respectively. Thus the Barnes--Wall lattices give the best known estimates in dimensions~$4$, $8$, and~$16$. However, starting from dimension~$32$, denser packings are known, for example, $Q_{32}$~\cite{Quebbemann} and $\mathrm{Ne}_{64}$~\cite{Nebe1998, Nebe1999}. Detailed information on the available results in various dimensions can be found in the monograph of Conway and Sloane~\cite[Introduction, pp.~xviii--xx, Table~I.1; Chapter~1, pp.~12--17, Tables~1.1--1.3]{ConwaySloane}.

One of the most elegant definitions of the Barnes--Wall lattices was given by Nebe, Rains, and Sloane~\cite{Nebe2002}. It is based on the so-called balanced Barnes--Wall lattices over the ring $\Z[\sqrt{2}] = \{a + \sqrt{2}b \mid a,b \in \Z\}$, whose generator matrices are constructed inductively by means of the tensor product of matrices. A similar construction, but over the Gaussian integers~$\Z[i]$, was introduced in the paper by Micciancio and Nicolosi~\cite{Micciancio}, where a family of bounded distance decoding algorithms for complex Barnes--Wall lattices was described.
The definition given in Forney's lecture notes~\cite[Chapter~14]{ForneyPrinciples} can be obtained by realification of the previous construction and becomes equivalent to the usual definition (namely, to~$2L_n$ in the notation of the paper~\cite{Nebe2002}) after applying the rotation $R = 
\left(
\begin{smallmatrix}
    1 & 1\\
    1 & -1
\end{smallmatrix}
\right)$ to each $2 \times 2$ block of the generator matrix.

Forney~\cite{ForneyCosetII} used the squaring construction for set partitions to obtain the Barnes--Wall lattices and their principal sublattices. Such a construction allows one to write down an explicit formula for their elements in terms of Reed--Muller codes, and it also clarifies why tensor products of matrices appear naturally in their generator matrices. Note that both real and complex Barnes--Wall lattices are studied in this paper, and in the complex case the formula via Reed--Muller codes takes a particularly simple form. Forney also considered the cubing construction which provides one of the many currently known constructions for the Golay code and the Leech lattice, see also~\cite{ConwaySloane1982}.

The construction of the Barnes--Wall lattices by means of the iterated doubling process suggests their connection with the Cayley--Dickson algebras. Throughout the introduction, by the Cayley--Dickson algebras we mean the so-called real Cayley--Dickson algebras of the main sequence~$\CD_n$, $n \geq 0$. They form a sequence of $2^n$-dimensional algebras which are defined inductively: $\CD_0 = \mathbb{R}$, and at each step the algebra~$\CD_{n+1}$ is obtained from~$\CD_n$ by applying the Cayley--Dickson doubling process with the parameter~$-1$. 
For $n = 3$ the algebra~$\CD_n$ becomes nonassociative, and for $n \geq 4$ it is not even alternative, that is, the identities $a(ab) = (aa)b$ and $(ba)a = b(aa)$ do not hold. The most widely known are the first members of this sequence --- for $n = 1$, $2$, $3$, and~$4$ these are the complex numbers, the quaternions, the octonions, and the sedenions, respectively.

One of the most interesting and natural questions concerning lattices embedded into the Cayley--Dickson algebras is whether they are closed under multiplication. It can be easily shown that in any Cayley--Dickson algebra the lattice of all elements with integer coordinates is closed, that is, it forms a subring. In the case of the quaternions such elements are called the Lipschitz quaternions. However, this is not the densest lattice with this property: the Hurwitz quaternion integers, whose coordinates are either all integers or all half-integers, are also closed under multiplication. They form a lattice which is geometrically equivalent to the $D_4$ lattice.

In the case of the octonions the situation turns out to be slightly more complicated. The Kirmse integers, whose construction is similar to that of the Hurwitz integers, turn out not to be closed under multiplication, but this can be fixed by swapping the coordinate at $1$ with any other coordinate for all elements~\cite[Section~9.2]{ConwaySmith}. This gives~$7$ copies of lattices which are closed under multiplication; each of them is geometrically equivalent to~$E_8$ and is called the Coxeter--Dickson integral octonions.

Among other important results devoted to the connection of the Cayley--Dickson algebras with various lattices, we mention Wilson's construction of the Leech lattice by means of the Coxeter--Dickson octonions~\cite{Wilson}. Besides, Dixon~\cite{Dixon1994Book,Dixon2011Book} explored the representation of the Leech lattice via $X$-products on the octonion algebra, as well as applications of this topic in physics.
Nevertheless, the lattices within the Cayley--Dickson algebras of dimension at least~$16$ have hardly been studied. For instance, Krau\ss{}har~\cite{Krausshar} gave a general definition of lattices with Cayley--Dickson multiplication, but all the main results were obtained for the octonion case. 

The aim of the present paper is to answer the question whether the Barnes--Wall lattices embedded into the Cayley--Dickson algebras of an appropriate dimension are closed under multiplication, and also to study the lattices formed by pairwise products of their elements. 

The structure of this paper is as follows: In Section~\ref{section:definitions} we describe the Cayley--Dickson doubling process and present some properties of the real Cayley--Dickson algebras of the main sequence. Section~\ref{section:Barnes-Wall} contains known and basic facts on the Barnes--Wall lattices: their construction, low-dimensional examples, and proofs of their main properties. We take the definition from the paper~\cite{Nebe2002} as a basis, with the only difference that we use an additional normalizing factor~$1/2$ in the formulas for the generator matrices. For each $n \in \N$ we first construct the balanced Barnes--Wall lattice~$M_n$ of dimension~$2^n$ over $\Z[\sqrt{2}]$, while its rational and purely irrational parts correspond to the ordinary Barnes--Wall lattices~$L_n$ and~$L'_n$ over~$\Z$, which satisfy the chain of inclusions $2L'_n \subset L_n \subset L'_n$.

In Theorems~\ref{theorem:balanced-closed} and~\ref{theorem:balanced-generated} of Section~\ref{section:balanced} we show that for $n \geq 2$ the balanced Barnes--Wall lattice~$M_n$ is closed under multiplication in~$\CD_n$, and all pairwise products of elements of~$M_n$ generate~$M_n$.

In Section~\ref{section:integer} we consider the lattices generated by the sets of pairwise products of elements of the ordinary Barnes--Wall lattices~$L_n$ and~$L'_n$. The main results of this section are Theorems~\ref{theorem:classic-closed}, \ref{theorem:classic-generated}, and~\ref{theorem:second-generated} which state that for all $n \geq 2$ we have
    \[
    L_n \cdot L_n = L_n, \qquad L_n \cdot L'_n = L'_n, \qquad L'_n \cdot L_n \subseteq L'_n, \qquad L'_n \cdot L'_n \subseteq L_n/2.
    \]
In particular, this means that~$L_n$ is closed under multiplication in~$\CD_n$, and all pairwise products of elements of~$L_n$ generate~$L_n$. Remarks~\ref{remark:dimension-2} and~\ref{remark:strict-inclusion} show that, although for $n = 2$ both inclusions above turn into equalities, for $n \geq 3$ they may be strict. Another interesting result of this section is Proposition~\ref{proposition:connection} which shows that~$L_n$ and~$L'_n$ are related to each other by multiplication by a certain nonzero element of~$\CD_n$.

\section{Cayley--Dickson algebras} \label{section:definitions}

In this section we give the definition and the main properties of the so-called real Cayley--Dickson algebras of the main sequence. This is the best-known sequence of the Cayley--Dickson algebras, whose first members are the real and the complex numbers, the quaternions, and the octonions. Note, however, that Cayley--Dickson algebras can be defined over any field, using arbitrary nonzero parameters at each doubling step, see the classical paper by Schafer~\cite{Schafer}. Among the authors who studied real Cayley--Dickson algebras of the main sequence, in particular, their zero divisors and alternative elements, we mention Moreno~\cite{Moreno, Moreno_alternative, Moreno_constructing} and Biss, Dugger, and Isaksen~\cite{Biss, Biss2}.

\begin{definition}
    We define the real algebras of the main sequence~$\CD_n$, $n \geq 0$, inductively:
    \begin{enumerate}
        \item $\CD_0 = \mathbb{R}$ is the field of real numbers with the identity involution $\bar{a} = a$ and the only basis element $e^{(0)}_0 = 1$;
        \item $\CD_{n+1}$ is the set of ordered pairs of elements of~$\CD_n$ with the component-wise addition, multiplication defined by
        \[
        (a,b)(c,d) =(ac - \bar{d}b,da+b\bar{c}),
        \]
        and the involution
        $$
        \qquad \overline{(a,b)}=(\bar{a},-b), \qquad a,b,c,d\in \CD_n.
        $$
        The standard basis elements $e^{(n+1)}_0, \dots, e^{(n+1)}_{2^{n+1}-1}$ of the algebra $\CD_{n+1}$ have the form
        \begin{equation*}
        e^{(n+1)}_m =
        \begin{cases}
        (e^{(n)}_m,0), & 0 \leq m \leq 2^n-1,\\
        (0,e^{(n)}_{m-2^n}), & 2^n \leq m \leq 2^{n+1}-1.
        \end{cases}
        \end{equation*}
    \end{enumerate}
\end{definition}

\begin{remark}[{\cite[p.~435]{Schafer}}]
For every $n \geq 0$, $\CD_n$ is a $2^n$-dimensional algebra over $\mathbb{R}$ with the unity $1 = e^{(n)}_0$.
Moreover, the involution on $\CD_n$ is regular, that is, $a + \bar{a} \in \mathbb{R}1$ and $a\bar{a} = \bar{a}a \in \mathbb{R}1$ for all $a \in \CD_n$.
\end{remark}

\begin{definition} \label{definition:real-imaginary-part}
For an arbitrary element $a \in \CD_n$ we define its {\em trace} as $\t(a) = a + \bar{a} \in \mathbb{R}$, and its {\em norm} as $\n(a) = a \bar{a} = \bar{a}a \in \mathbb{R}$.
We also define the symmetric bilinear form corresponding to the quadratic form $\n(\cdot)$ by the formula
\[
\n(a, b) = \n(a+b) - \n(a) - \n(b), \quad a,b \in \CD_n.
\]
\end{definition}

\begin{proposition}[{\cite[p.~435]{Schafer}}] \label{proposition:real-cayley-dickson-properties}
The trace and the norm of an element $(a,b) \in \CD_{n+1}$ can be computed inductively by using the following relations:
\begin{align*}
	\t((a,b)) &= \t(a),\\
	\n((a,b)) &= \n(a) + \n(b).
\end{align*}
\end{proposition}

\begin{corollary}[{\cite[Proposition~3.2]{Biss}}]
The norm $\n(\cdot)$ on $\CD_n$ equals the square of the Euclidean norm in the standard basis, and the bilinear form $\n(\cdot,\cdot)$ equals twice the Euclidean inner product.
\end{corollary}

\begin{remark}
The norm on the real algebras of the main sequence is often defined as $\sqrt{a\bar{a}}$, in contrast to $\n(a) = a\bar{a}$ used in the present paper. Besides, instead of the trace $\t(a)$ and the bilinear form $\n(a,b)$, one often considers the real part $\mathrm{Re}(a) = \t(a)/2$ and the Euclidean inner product $\n(a,b)/2$. 

However, we prefer the notations $\n(a)$ and $\n(a,b)$, since they are universal and can be applied to Cayley--Dickson algebras over arbitrary fields, including fields of characteristic~$2$. They are also standard in the theory of composition algebras. Moreover, the use of precisely these notations will significantly shorten the computations in the proofs of the main results of the present paper.
\end{remark}

\begin{proposition}[{\cite[p.~440]{Schafer}}] \label{proposition:lambda-form}
For all $a, b \in \CD_n$ we have
\begin{enumerate}
    \item $\n(a, a) = 2\n(a)$,
    \item $\n(a, 1) = \t(a)$,
    \item $\n(a,b) = a \bar{b} + b \bar{a} = \bar{a} b + \bar{b} a = \t(a\bar{b})$,
    \item $\n(a, b) =\n(\bar{a}, \bar{b})$.
\end{enumerate}
\end{proposition}

It is well known that the Cayley--Dickson algebras are commutative only for $n \leq 1$, associative for $n \leq 2$, alternative (hence, composition algebras) for $n \leq 3$, and for $n \geq 4$ they are neither commutative, nor alternative, nor composition algebras, see~\cite[p.~10]{Baez}. Due to the loss of alternativity and the composition property, zero divisors appear in the Cayley--Dickson algebras starting from dimension~$16$. Nevertheless, all Cayley--Dickson algebras are flexible, that is, satisfy the identity $(ab)a = a(ba)$, see~\cite[p.~437, Theorem~1]{Schafer}.

\begin{example}
The algebras of the complex numbers ($\mathbb{C}$), the quaternions ($\mathbb{H}$), the octonions ($\mathbb{O}$), and the sedenions ($\mathbb{S}$) are the real algebras of the main sequence for $n = 1$, $2$, $3$, and $4$, respectively, see~\cite{Baez}.

The octonion algebra $\mathbb{O} = \CD_3$ is an eight-dimensional algebra over~$\mathbb{R}$ with the standard basis $1,e_1,\dots,e_7$. The involution on $\mathbb{O}$ is given by the formula
$$
\overline{a_0 + a_1e_1 + \dots + a_7e_7} = a_0 - a_1e_1 - \dots - a_7e_7,
$$
and the multiplication is given by Table~\ref{table:octonions-mult}.

\begin{table}[ht]
\centering
$
\begin{array}{|c|cccccccc|}
\hline
\times&1&e_1&e_2&e_3&e_4&e_5&e_6&e_7\\\hline
1&1&e_1&e_2&e_3&e_4&e_5&e_6&e_7\\
e_1&e_1&-1& e_3& -e_2& e_5& -e_4& -e_7& e_6\\
e_2&e_2&-e_3& -1& e_1& e_6& e_7& -e_4& -e_5\\
e_3&e_3&e_2& -e_1& -1& e_7& -e_6& e_5& -e_4\\
e_4&e_4&-e_5& -e_6& -e_7& -1& e_1& e_2& e_3\\
e_5&e_5&e_4& -e_7& e_6& -e_1& -1& -e_3& e_2\\
e_6&e_6&e_7& e_4& -e_5& -e_2& e_3& -1& -e_1\\
e_7&e_7&-e_6& e_5& e_4& -e_3& -e_2& e_1& -1\\\hline
\end{array}
$
\caption{\label{table:octonions-mult} The multiplication table of the basis octonions.}
\end{table}
\end{example}

Since $n =3$, $\mathbb{O}$ is a noncommutative, nonassociative, but alternative algebra without zero divisors. The basis $1, e_1, \dots, e_7$ is orthonormal with respect to the Euclidean inner product $\n(\cdot,\cdot)/2$.

\section{Barnes--Wall lattices} \label{section:Barnes-Wall}

In geometry, a {\em lattice} is a discrete additive subgroup of $\mathbb{R}^m$, $m \in \N$. It is isomorphic to a free abelian group of rank $n \leq m$, and in this case the number $n$ is called the {\em dimension} of the lattice. A {\em basis} of a lattice is its basis as a free abelian group, and all elements of the basis are linearly independent over~$\mathbb{R}$. A lattice is called {\em full} if $n = m$, that is, the basis of the lattice is a basis of $\mathbb{R}^m$. A {\em generator matrix} of a lattice is a matrix whose rows form its basis. Other basic notions related to lattices can be found in~\cite[Chapter~1, pp.~4--5]{ConwaySloane}. 

The standard definition allows one to define lattices over $\Z$, but it can be generalized to the case of an arbitrary domain $R$, see, e.g.,~\cite[Chapter~9]{Voight}. Namely, consider its field of fractions $F$ and a finite-dimensional vector space $V$ over $F$. Then a lattice~$M$ over $R$ is a finitely generated $R$-submodule $M \subset V$. A lattice $M$ is called full if $V = FM$.

We follow the approach of~\cite{Nebe2002} for the construction of the Barnes--Wall lattices. Put
\[
G = \begin{pmatrix}
    \sqrt{2} & 0\\
    1 & 1
\end{pmatrix}
\quad \text{and} \quad
G_n = \frac{1}{2}G^{\otimes n}, \quad n \geq 0,
\]
where~$\otimes$ is the tensor (Kronecker) product of matrices. Here we assume that $G^{\otimes 0} = \big( 1 \big)$. We emphasize that we use an additional normalizing factor~$1/2$ in the formula for~$G_n$ in order to simplify the statements of the main results of the present paper.

\begin{definition} \label{definition:Barnes-Wall-lattices}
    The {\em balanced Barnes--Wall lattice} $M_n$ is a $2^n$-dimensional lattice over $\Z[\sqrt{2}] = \Z + \sqrt{2}\Z$ whose generator matrix is~$G_n$. The rational part of~$M_n$ is the lattice~$L_n$, and the purely irrational part is the lattice $\sqrt{2} L'_n$. The lattices~$L_n$ and~$L'_n$ are the {\em Barnes--Wall lattices} of dimension~$2^n$ over~$\Z$.
\end{definition}

Definition~\ref{definition:Barnes-Wall-lattices} yields
\[
L_n + \sqrt{2} L'_n = M_n \supset \sqrt{2}M_n = 2L'_n + \sqrt{2} L_n.
\]
This immediately implies the well-known chain of inclusions
\begin{equation} \label{equation:inclusion}
    2L'_n \subset L_n \subset L'_n.
\end{equation}
It follows from Proposition~\ref{proposition:common-types} below that $L_n$ is a sublattice of index $2^{2^{n-1}}$ in $L'_n$, and $2 L'_n$ is a sublattice of index $2^{2^{n-1}}$ in $L_n$. Extending the chain~\eqref{equation:inclusion} to the left, we obtain
\[
\ldots \subset 2^k L_n \subset 2^k L'_n \subset \ldots \subset 2L_n \subset 2L'_n \subset L_n \subset L'_n.
\]
A similar chain was obtained in Forney's paper~\cite[pp.~1167--1169]{ForneyCosetII} by means of the squaring construction for set partitions. Here the lattice~$2L'_n$ was replaced by the lattice~$RL_n$ which is obtained from~$L_n$ by applying the rotation $R = 
\left(
\begin{smallmatrix}
    1 & 1\\
    1 & -1
\end{smallmatrix}
\right)$ to each $2 \times 2$ block of the generator matrix in a suitable basis.

\begin{notation}
    We denote by~$H_n$ and~$H'_n$ the generator matrices of the lattices~$L_n$ and~$L'_n$, respectively, whose rows are obtained by multiplication of the corresponding rows of the matrix~$G_n$ by positive numbers.
\end{notation}

The following example demonstrates the structure of the matrices~$H_n$ and~$H'_n$ for $n \leq 3$, and Propositions~\ref{proposition:common-types} and~\ref{proposition:connection-two-lattices} fully describe their rows in an arbitrary dimension.

\begin{example} \label{example:low-dimension}
    For $n = 0$ the lattices~$M_0$, $L_0$, and~$L'_0$ have the common generator matrix
    \[
    G_0 = H_0 = H'_0 = \begin{pmatrix}
    \frac{1}{2}
    \end{pmatrix}.
    \]
    For $n = 1$ the generator matrices of~$M_1$, $L_1$, and~$L'_1$ are 
    \[
    G_1 = \begin{pmatrix}
        \frac{\sqrt{2}}{2} & 0\\
        \frac{1}{2} & \frac{1}{2}
    \end{pmatrix},
    \quad
    H_1 = \begin{pmatrix}
        1 & 0\\
        \frac{1}{2} & \frac{1}{2}
    \end{pmatrix}
    \quad \text{and} \quad
    H'_1 = \begin{pmatrix}
        \frac{1}{2} & 0\\
        \frac{1}{2} & \frac{1}{2}
    \end{pmatrix}
    \]
    respectively. For $n = 2$ the generator matrices are
    \[
    G_2 = \begin{pmatrix}
        1 & 0 & 0 & 0\\
        \frac{\sqrt{2}}{2} & \frac{\sqrt{2}}{2} & 0 & 0\\
        \frac{\sqrt{2}}{2} & 0 & \frac{\sqrt{2}}{2} & 0\\
        \frac{1}{2} & \frac{1}{2} & \frac{1}{2} & \frac{1}{2}
    \end{pmatrix},
    \quad
    H_2 = \begin{pmatrix}
        1 & 0 & 0 & 0\\
        1 & 1 & 0 & 0\\
        1 & 0 & 1 & 0\\
        \frac{1}{2} & \frac{1}{2} & \frac{1}{2} & \frac{1}{2}
    \end{pmatrix}
    \quad \text{and} \quad
    H'_2 = \begin{pmatrix}
        1 & 0 & 0 & 0\\
        \frac{1}{2} & \frac{1}{2} & 0 & 0\\
        \frac{1}{2} & 0 & \frac{1}{2} & 0\\
        \frac{1}{2} & \frac{1}{2} & \frac{1}{2} & \frac{1}{2}
    \end{pmatrix}.
    \]
    Therefore, $L_2$ and~$L'_2$ are geometrically equivalent to the lattice~$D_4$, see~\cite[Chapter~4, Eqs. (86) and (90)]{ConwaySloane}. Next, we write down the generator matrices for $n = 3$ and the lower triangular Hermite normal forms for~$H_3$ and~$H'_3$ (computed row-wise):
    \begin{gather*}
        G_3 = \left(\begin{array}{cccccccc}
        \sqrt{2} & 0 & 0 & 0 & 0 & 0 & 0 & 0\\
        1 & 1 & 0 & 0 & 0 & 0 & 0 & 0\\
        1 & 0 & 1 & 0 & 0 & 0 & 0 & 0\\
        \frac{\sqrt{2}}{2} & \frac{\sqrt{2}}{2} & \frac{\sqrt{2}}{2} & \frac{\sqrt{2}}{2} & 0 & 0 & 0 & 0\\
        1 & 0 & 0 & 0 & 1 & 0 & 0 & 0\\
        \frac{\sqrt{2}}{2} & \frac{\sqrt{2}}{2} & 0 & 0 & \frac{\sqrt{2}}{2} & \frac{\sqrt{2}}{2} & 0 & 0\\
        \frac{\sqrt{2}}{2} & 0 & \frac{\sqrt{2}}{2} & 0 & \frac{\sqrt{2}}{2} & 0 & \frac{\sqrt{2}}{2} & 0\\
        \frac{1}{2} & \frac{1}{2} & \frac{1}{2} & \frac{1}{2} & \frac{1}{2} & \frac{1}{2} & \frac{1}{2} & \frac{1}{2}\\
        \end{array}\right),\\
        \begin{aligned}
        H_3 = \left(\begin{array}{cccccccc}
        2 & 0 & 0 & 0 & 0 & 0 & 0 & 0 \\
        1 & 1 & 0 & 0 & 0 & 0 & 0 & 0 \\
        1 & 0 & 1 & 0 & 0 & 0 & 0 & 0 \\
        1 & 1 & 1 & 1 & 0 & 0 & 0 & 0 \\
        1 & 0 & 0 & 0 & 1 & 0 & 0 & 0\\
        1 & 1 & 0 & 0 & 1 & 1 & 0 & 0\\
        1 & 0 & 1 & 0 & 1 & 0 & 1 & 0\\
        \frac{1}{2} & \frac{1}{2} & \frac{1}{2} & \frac{1}{2} & \frac{1}{2} & \frac{1}{2} & \frac{1}{2} & \frac{1}{2} \\
        \end{array}\right)
        &\to \widetilde{H}_3 = \left(\begin{array}{cccccccc}
        2 & 0 & 0 & 0 & 0 & 0 & 0 & 0 \\
        1 & 1 & 0 & 0 & 0 & 0 & 0 & 0 \\
        1 & 0 & 1 & 0 & 0 & 0 & 0 & 0 \\
        1 & 0 & 0 & 1 & 0 & 0 & 0 & 0 \\
        1 & 0 & 0 & 0 & 1 & 0 & 0 & 0\\
        1 & 0 & 0 & 0 & 0 & 1 & 0 & 0\\
        1 & 0 & 0 & 0 & 0 & 0 & 1 & 0\\
        \frac{1}{2} & \frac{1}{2} & \frac{1}{2} & \frac{1}{2} & \frac{1}{2} & \frac{1}{2} & \frac{1}{2} & \frac{1}{2} \\
        \end{array}\right),\\
        H'_3 = \left(\begin{array}{cccccccc}
        1 & 0 & 0 & 0 & 0 & 0 & 0 & 0 \\
        1 & 1 & 0 & 0 & 0 & 0 & 0 & 0 \\
        1 & 0 & 1 & 0 & 0 & 0 & 0 & 0 \\
        \frac{1}{2} & \frac{1}{2} & \frac{1}{2} & \frac{1}{2} & 0 & 0 & 0 & 0 \\
        1 & 0 & 0 & 0 & 1 & 0 & 0 & 0\\
        \frac{1}{2} & \frac{1}{2} & 0 & 0 & \frac{1}{2} & \frac{1}{2} & 0 & 0\\
        \frac{1}{2} & 0 & \frac{1}{2} & 0 & \frac{1}{2} & 0 & \frac{1}{2} & 0\\
        \frac{1}{2} & \frac{1}{2} & \frac{1}{2} & \frac{1}{2} & \frac{1}{2} & \frac{1}{2} & \frac{1}{2} & \frac{1}{2} \\
        \end{array}\right)
        &\to \widetilde{H}'_3 = \left(\begin{array}{cccccccc}
        1 & 0 & 0 & 0 & 0 & 0 & 0 & 0 \\
        0 & 1 & 0 & 0 & 0 & 0 & 0 & 0 \\
        0 & 0 & 1 & 0 & 0 & 0 & 0 & 0 \\
        \frac{1}{2} & \frac{1}{2} & \frac{1}{2} & \frac{1}{2} & 0 & 0 & 0 & 0 \\
        0 & 0 & 0 & 0 & 1 & 0 & 0 & 0\\
        \frac{1}{2} & \frac{1}{2} & 0 & 0 & \frac{1}{2} & \frac{1}{2} & 0 & 0\\
        \frac{1}{2} & 0 & \frac{1}{2} & 0 & \frac{1}{2} & 0 & \frac{1}{2} & 0\\
        0 & \frac{1}{2} & \frac{1}{2} & 0 & \frac{1}{2} & 0 & 0 & \frac{1}{2} \\
        \end{array}\right).\\
        \end{aligned}
    \end{gather*}

    Due to the additional factor $1/2$, the matrix~$\widetilde{H}_3$ is a generator matrix of the lattice~$E_8$ written in the standard form via the even coordinate system: all coordinates of an arbitrary element are either integers or half-integers simultaneously, and their sum is even.
\end{example}

Relying on the previous example, it is not hard to describe the properties of the generator matrices of the Barnes--Wall lattices in an arbitrary dimension.

\begin{proposition} \label{proposition:main-properties}
    We list the main properties of the matrix $G_n$:
    \begin{enumerate}
        \item $G_n$ is a lower triangular matrix with nonzero diagonal entries;
        \item in each row of $G_n$ all entries belong simultaneously either to $\Z/2$ or to $\sqrt{2}\Z/2$;
        \item the last row of $G_n$ equals $(1/2, \dots, 1/2)$;
        \item the first row of $G_n$ equals $((\sqrt{2})^{n-2}, 0, \dots, 0)$.
    \end{enumerate}
\end{proposition}

\begin{proof}
    This is easily obtained by induction on~$n$.
\end{proof}

\begin{remark} \label{remark:rational-part}
    It follows from Definition~\ref{definition:Barnes-Wall-lattices} that $L_n, \sqrt{2} L'_n \subset M_n$.  The coordinates of an arbitrary element of the lattice $M_n$ belong to $\Z[\sqrt{2}]/2$, hence the coordinates of all elements of $L_n$ and $L'_n$ belong to $\Z/2$. The rational and the purely irrational parts of the lattice~$M_n$ are given by the mappings
    \[
        a + \sqrt{2}b \mapsto a + 2b \quad \text{and} \quad a + \sqrt{2}b \mapsto \sqrt{2}(a + b), \quad a,b \in \Z/2,
    \]
    respectively, applied coordinate-wise.
\end{remark}

\begin{definition}
    We say that a row of~$G_n$ has {\em type I} if all its coordinates belong to~$\Z/2$, and {\em type II} if all its coordinates belong to $\sqrt{2}\Z/2$. The type of an arbitrary row of~$H_n$ or~$H'_n$ coincides with the type of the corresponding row of~$G_n$.
\end{definition}

\begin{proposition} \label{proposition:types}
    \leavevmode
    \begin{enumerate}
        \item The row~$(a,a)$ has the same type as~$a$.
        \item The rows~$a$ and $(\sqrt{2}a,0)$ have different types.
        \item The last row always has type I.
        \item The first row has type II if and only if~$n$ is odd.
    \end{enumerate}
\end{proposition}

\begin{proof}
    Items (1--2) follow immediately from the definition of types, and items (3--4) follow from Proposition~\ref{proposition:main-properties}(3--4).
\end{proof}

\begin{proposition} \label{proposition:common-types}
    The rows of type~I are common for the generator matrices~$G_n$, $H_n$, and~$H'_n$. The rows of type II in~$H_n$ are~$\sqrt{2}$ times greater than the corresponding rows of~$G_n$ and twice as great as the corresponding rows of~$H'_n$.
\end{proposition}

\begin{proof}
    If a row~$a$ of the matrix~$G_n$ has type~I, that is, $a \in (\Z/2)^{2^n}$, then it corresponds to~$a$ and~$\sqrt{2}a$ in the rational and the purely irrational parts of~$M_n$, respectively. Hence in~$H_n$ and~$H'_n$ the row~$a$ remains unchanged. If~$a$ has type~II, then $a = \sqrt{2}b$, where $b \in (\Z/2)^{2^n}$, so it corresponds to~$2b$ and~$\sqrt{2}b$ in the rational and the purely irrational parts of~$M_n$. Consequently, the rows of~$H_n$ and~$H'_n$ are~$2b$ and~$b$, respectively.
\end{proof}

The following statements show how to construct the elements of the $2^{n+1}$-dimensional Barnes--Wall lattices from the elements of the $2^n$-dimensional lattices.

\begin{proposition} \label{proposition:new-elements}
    If $x \in M_n$, then $(\sqrt{2}x,0)$, $(0,\sqrt{2}x)$, $(x,x)$, and~$(x,-x)$ belong to~$M_{n+1}$.
\end{proposition}

\begin{proof}
    By linearity, it suffices to prove this statement in the case when~$x$ is a row of the matrix~$G_n$. According to the definition,
    \[
    G_{n+1} = \begin{pmatrix}
        \sqrt{2} G_n & 0\\
        G_n & G_n
    \end{pmatrix}.
    \]
    Consequently, an arbitrary row of~$G_{n+1}$ has the form $(\sqrt{2}a,0)$ or~$(a,a)$, where~$a$ is some row of~$G_n$. Hence $(\sqrt{2}x,0), (x,x) \in M_{n+1}$. Besides, 
    \begin{align*}
        (0,\sqrt{2}x) &= \sqrt{2}(x,x) - (\sqrt{2}x,0) \in M_{n+1},\\
        (x,-x) &= (x,x) - \sqrt{2}(0,\sqrt{2}x) \in M_{n+1}. \tag*{\qedhere}
    \end{align*}
\end{proof}

\begin{proposition} \label{proposition:connection-two-lattices}
    The rows of the generator matrix $H_{n+1}$ of the lattice $L_{n+1}$ have the form $(x,x)$ and $(2x',0)$, and the rows of the generator matrix $H'_{n+1}$ of the lattice $L'_{n+1}$ have the form $(x',x')$ and $(x,0)$, where $x$ is an arbitrary row of $H_n$, and $x'$ is an arbitrary row of $H'_n$.
\end{proposition}

\begin{proof}
    An arbitrary row $a$ of the matrix $G_n$ corresponds to the rows $(a,a)$ and $(\sqrt{2}a,0)$ of the matrix $G_{n+1}$.
    Suppose first that $a$ has type~I, then $(a,a)$ also has type~I, and $(\sqrt{2}a,0)$ has type~II.
    By Proposition~\ref{proposition:common-types}, the row $a$ is common for $H_n$ and $H'_n$, the row $(a,a)$ is common for $H_{n+1}$ and $H'_{n+1}$, and the row $(\sqrt{2}a,0)$ corresponds to $(2a,0)$ in $H_{n+1}$ and to $(a,0)$ in $H'_{n+1}$. 

    Now let $a$ have type~II, then $(a,a)$ has type~II, and $(\sqrt{2}a,0)$ has type~I.
    The row~$a$ corresponds to $\sqrt{2}a$ in $H_n$ and to $a/\sqrt{2}$ in $H'_n$, the row $(a,a)$ corresponds to $\sqrt{2}(a,a)$ in $H_{n+1}$ and to $(a,a)/\sqrt{2}$ in $H'_{n+1}$, and the row $(\sqrt{2}a,0)$ is common for $H_{n+1}$ and $H'_{n+1}$. 
    
    Combining these two cases, we obtain the required statement.
\end{proof}

\begin{corollary} \label{corollary:new-elements}
    Let $x \in L_n$, $x' \in L'_n$ be arbitrary elements. Then
    \begin{align*}
        (x,x), (x,-x), (2x',0), (0,2x') &\in L_{n+1},\\
        (x',x'), (x',-x'), (x,0), (0,x) &\in L'_{n+1}.
    \end{align*}
\end{corollary}

\begin{proof}
    It follows from the chain of inclusions~\eqref{equation:inclusion} that $2x' \in 2L'_n \subset L_n$ and $x \in L_n \subset L'_n$. By Proposition~\ref{proposition:connection-two-lattices}, it suffices to note that
    \begin{align*}
        &(x,-x) = (x,x) - (0,2x) \in L_{n+1}, &
        &(0,2x') = (2x',2x') - (2x',0) \in L_{n+1},\\
        &(x',-x') = (x',x') - (0,2x') \in L'_{n+1}, &
        &(0,x) = (x,x) - (x,0) \in L'_{n+1}. \tag*{\qedhere}
    \end{align*}
\end{proof}

\begin{corollary} \label{corollary:new-elements-Ln}
    If $x \in L_n$, then $(2x,0), (0,2x) \in L_{n+1}$.
\end{corollary}

\begin{proof}
    This follows immediately from Corollary~\ref{corollary:new-elements} and the inclusion $L_n \subset L'_n$.
\end{proof}

\begin{corollary} \label{corollary:new-elements-types}
    Each row $a$ of the generator matrix $H_n$ corresponds to two rows of the matrix $H_{n+1}$: these are $(a,a)$ and either $(a,0)$ or $(2a,0)$. In the first case $a$ has type II, while in the second case it has type I.
\end{corollary}

\begin{proof}
    Follows immediately from Propositions~\ref{proposition:common-types} and~\ref{proposition:connection-two-lattices}.
\end{proof}

\section{Products of elements of the balanced Barnes--Wall lattice} \label{section:balanced}

Consider the embedding of the lattice $M_n$ into the Cayley--Dickson algebra of the main sequence $\CD_n$ such that the elements of the lattice correspond to the coordinates of the algebra elements in the standard basis $1 = e_0, e_1, \dots, e_{2^n-1}$. In this section we show that for $n \geq 2$ the lattice $M_n$ is closed under multiplication, and, moreover, all pairwise products of elements of $M_n$ generate $M_n$ over $\Z[\sqrt{2}]$.

\begin{remark} \label{remark:not-closed-dimension-one}
    The lattices $M_1$ and $L_1$ are not closed under multiplication, since the element $z = (1+i)/2$, which corresponds to the last row of the matrix $G_1$, satisfies $z^2 = i/2 \notin M_1$. The lattice $L'_1$ is not closed under multiplication either, since the element $y = 1/2$, which corresponds to the first row of the matrix $H'_1$, has $y^2 = 1/4 \notin L'_1$.
\end{remark}

\begin{remark} \label{remark:closed-dimension-two}
    One can directly verify that the lattice $M_2$ is closed under multiplication by writing down the multiplication table of the rows of its generator matrix $G_2$, see Table~\ref{table:M2-multiplication}.
    Moreover, $1 \in M_2$, which immediately implies that all pairwise products of elements of $M_2$ generate~$M_2$.
    \begin{table}[H]
    \centering
    $
    \begin{array}{|c||c|c|c|c|}
    \hline
    \vphantom{\Bigg|}\times&1&\dfrac{1+i}{\sqrt{2}}&\dfrac{1+j}{\sqrt{2}}&\dfrac{1+i+j+k}{2}\\\hline\hline
    \vphantom{\Bigg|}1&1&\dfrac{1+i}{\sqrt{2}}&\dfrac{1+j}{\sqrt{2}}&\dfrac{1+i+j+k}{2}\\\hline
    \vphantom{\Bigg|}\dfrac{1+i}{\sqrt{2}}&\dfrac{1+i}{\sqrt{2}}&i&\dfrac{1+i+j+k}{2}&\dfrac{i+k}{\sqrt{2}}\\\hline
    \vphantom{\Bigg|}\dfrac{1+j}{\sqrt{2}}&\dfrac{1+j}{\sqrt{2}}&\dfrac{1+i+j-k}{2}&j&\dfrac{i+j}{\sqrt{2}}\\\hline
    \vphantom{\Bigg|}\dfrac{1+i+j+k}{2}&\dfrac{1+i+j+k}{2}&\dfrac{i+j}{\sqrt{2}}&\dfrac{j+k}{\sqrt{2}}&\dfrac{-1+i+j+k}{2}\\\hline
    \end{array}
    $
    \caption{\label{table:M2-multiplication} The multiplication table of the rows of $G_2$.}
    \end{table}
\end{remark}

\begin{proposition} \label{proposition:rows-norm}
    Let $x$ and $y$ be rows of the generator matrix~$G_n$ of the lattice~$M_n$, $n \geq 0$. Then $\n(x,y) = (\sqrt{2})^m$ for some integer $m \geq n - 2$, and $m = n - 2$ if and only if the sum of the indices of these rows equals $2^n + 1$, that is, the rows $x$ and $y$ are located symmetrically with respect to the horizontal axis of the matrix $G_n$.
\end{proposition}

\begin{proof}
    We use induction on $n$. If $n=0$, then $x = y = (1/2)$, whence $\n(x,y) = 1/2 = (\sqrt{2})^{-2}$, since $\n(x,y)$ is twice the inner product of $x$ and $y$.
        
        Now suppose that the statement holds for some $n \geq 0$. It follows from the equalities
        \begin{align*}
            \n((\sqrt{2}x,0),(\sqrt{2}y,0)) &= 2\n(x,y),\\
            \n((\sqrt{2}x,0),(y,y)) &= \sqrt{2}\n(x,y),\\
            \n((x,x),(y,y)) &= 2\n(x,y)
        \end{align*}
        that the value of $\n(\cdot,\cdot)$ on pairs of rows of $G_{n+1}$ equals $(\sqrt{2})^{m'}$, where $m' \geq m + 1 \geq (n+1)-2$. The equality $m' = (n+1)-2$ is attained only in the second of these three cases, provided that $\n(x,y) = (\sqrt{2})^{n-2}$. The index of the row $(\sqrt{2}x,0)$ in $G_{n+1}$ coincides with the index of the row $x$ in $G_n$, and the index of the row $(y,y)$ in $G_{n+1}$ exceeds the index of $y$ in $G_n$ by $2^n$. Thus the equality is attained if and only if the sum of the row indices equals $(2^n + 1) + 2^n = 2^{n+1} + 1$.
\end{proof}

\begin{proposition} \label{proposition:trace}
    For all $x,y \in M_n$, $n \geq 0$, we have
    \begin{enumerate}
        \item $\t(x) \in \Z[\sqrt{2}]$,
        \item $\n(x,y)$ is divisible by $(\sqrt{2})^{n-2}$, that is, there exists $k \in \Z[\sqrt{2}]$ such that $\n(x,y) = k(\sqrt{2})^{n-2}$.
    \end{enumerate}
\end{proposition}

\begin{proof}
    \begin{enumerate}
        \item The trace $\t(x)$ equals twice the coordinate of $x$ at $1 = e_0$, which belongs to $\Z[\sqrt{2}]/2$ by Remark~\ref{remark:rational-part}.
        \item This follows from Proposition~\ref{proposition:rows-norm} and the bilinearity of the form $\n(\cdot,\cdot)$.  \qedhere
    \end{enumerate}
\end{proof}

\begin{theorem} \label{theorem:balanced-closed}
    If $n \geq 2$, then the lattice $M_n$ is closed under multiplication in $\CD_n$.
\end{theorem}

\begin{proof}
    We use induction on $n$. Remark~\ref{remark:closed-dimension-two} demonstrates that the lattice $M_2$ is closed under multiplication. Now suppose that the statement holds for $M_n$, and let us prove it for $M_{n+1}$. It suffices to show that all pairwise products of rows of $G_{n+1}$ belong to $M_{n+1}$. Consider four types of pairwise products, assuming that $a$ and $b$ are arbitrary rows of~$G_n$:
    \begin{itemize}
        \item $(\sqrt{2}a,0)(\sqrt{2}b,0) = (2ab,0)$. By the induction hypothesis, $ab \in M_n$, hence $(\sqrt{2}ab,0) \in M_{n+1}$, so the product under consideration belongs to $\sqrt{2} M_{n+1}$.
        \item $(b,b)(\sqrt{2}a,0) = \sqrt{2}(ba,b\bar{a}) = \sqrt{2}(ba,\t(a)b-ba) = \t(a)(0,\sqrt{2}b) + \sqrt{2}(ba,-ba)$. By the induction hypothesis, $ba \in M_n$, hence it follows from Proposition~\ref{proposition:new-elements} that $(0,\sqrt{2}b), (ba,-ba) \in M_{n+1}$. Since $\t(a) \in \Z[\sqrt{2}]$ by Proposition~\ref{proposition:trace}(1), the product under consideration belongs to $M_{n+1}$.
        \item $(a,a)(b,b) = (ab - \bar{b}a,ba+a\bar{b}) = (ab - \t(b)a + ba, ba + \t(b)a - ab) = (ba,ba) + (ab,-ab) - \t(b)(a,-a)$. Applying the induction hypothesis and Propositions~\ref{proposition:new-elements} and~\ref{proposition:trace}(1) again, we obtain the required statement.
        \item $(\sqrt{2}a,0)(b,b) = \sqrt{2}(ab,ba)$. Further, expressing $ab$ from the equality
        \[
        \n(a,b) = a \bar{b} + b \bar{a} = \t(b)a + \t(a)b - ab - ba,
        \]
        we obtain that
        \[
        (\sqrt{2}a,0)(b,b) = \t(b)(\sqrt{2}a,0) + \t(a)(\sqrt{2}b,0) - \sqrt{2}(ba,-ba) - \sqrt{2} \n(a,b) (1,0).
        \]
        Proposition~\ref{proposition:trace}(2) implies that $\n(a,b)$ is divisible by $(\sqrt{2})^{n-2}$. Since, by Proposition~\ref{proposition:main-properties}(4), the first row of~$G_{n+1}$ equals $((\sqrt{2})^{n-1}, 0, \dots, 0)$, we have $\sqrt{2} \n(a,b) (1,0) \in M_{n+1}$. The remaining summands are treated similarly to the previous cases.
        \qedhere
    \end{itemize}
\end{proof}

\begin{theorem} \label{theorem:balanced-generated}
    If $n \geq 2$, then all pairwise products of elements of~$M_n$ generate~$M_n$ over $\Z[\sqrt{2}]$.
\end{theorem}

\begin{proof}
    For convenience of notation, we prove the statement for~$M_{n+1}$, $n \in \N$. Denote by~$V_{n+1}$ the sublattice of~$M_{n+1}$ over $\Z[\sqrt{2}]$ generated by all pairwise products of elements of~$M_{n+1}$. It suffices to show that all rows of the generator matrix~$G_{n+1}$ belong to~$V_{n+1}$. For $a, b \in M_n$ consider the equality
    \begin{align*}
        V_{n+1} \ni (b,b)(0,\sqrt{2}a) + (a,a)(\sqrt{2}b,0) &= \sqrt{2}((-\bar{a}b,ab)+(ab,a\bar{b})) \\
        &= \sqrt{2} (\t(b)(0,a) - \t(a)(b,0) + (2ab,0)).
    \end{align*}
    Since $(2ab,0) = (\sqrt{2}a,0)(\sqrt{2}b,0) \in V_{n+1}$, we obtain that 
    \begin{equation} \label{equation:generated-sum}
        \t(b)(0,\sqrt{2}a) - \t(a)(\sqrt{2}b,0) \in V_{n+1}.
    \end{equation}

    Let $z = (1/2, \dots, 1/2)$ be the last row of $G_n$, see Proposition~\ref{proposition:main-properties}(3). Then $\t(z) = 1$, and we have the equality
    \begin{align*}
        V_{n+1} \ni (\sqrt{2}z,0)(z,z) - (z,z)(0,\sqrt{2}z) &= \sqrt{2}((z^2,z^2)-(-\bar{z}z,z^2)) \\
        &= \t(z)(\sqrt{2}z,0) = (\sqrt{2}z,0).
    \end{align*}
    Putting $a = b = z$ in Eq.~\eqref{equation:generated-sum}, we obtain that $(0,\sqrt{2}z) \in V_{n+1}$.

    Further, if $b = z$ and $a \in M_n$ is an arbitrary element, then $\t(b) = 1$, and Eq.~\eqref{equation:generated-sum} implies $(0,\sqrt{2}a) \in V_{n+1}$. Similarly, if $a = z$ and $b \in M_n$ is an arbitrary element, then $(\sqrt{2}b,0) \in V_{n+1}$.

    It remains to show that $V_{n+1}$ contains those rows of $G_{n+1}$ which have the form $(c,c)$. For arbitrary $a,b \in M_n$ consider the equality
    \begin{align*}
        V_{n+1} \ni (a,a)(b,b) &= (ab - \bar{b}a,ba + a\bar{b}) \\
        &= (\t(a)b - (\bar{a}b + \bar{b}a),-\t(a)b + 2ba + (b\bar{a}+a\bar{b})) \\
        &= (\t(a)b - \n(a,b),-\t(a)b + 2ba + \n(a,b)) \\
        &= \t(a)(b,b) + \n(a,b)(1,1) - 2\t(a)(0,b) + 2(0,ba) - 2(\n(a,b),0).
    \end{align*}
    By Proposition~\ref{proposition:trace}(2), $\n(a,b)$ is divisible by $(\sqrt{2})^{n-2}$, and, by Proposition~\ref{proposition:main-properties}(4), the first row of $G_{n}$ equals $((\sqrt{2})^{n-2}, 0, \dots, 0)$. Besides, it follows from Theorem~\ref{theorem:balanced-closed} that $ba \in M_n$. By what has been proved above, we have 
    \[
    (0,\sqrt{2}b), (0,\sqrt{2}ba), (\sqrt{2}\n(a,b),0) \in V_{n+1}.
    \]
    Proposition~\ref{proposition:trace}(1) implies that $\t(a) \in \Z[\sqrt{2}]$, hence
    \[
    \t(a)(b,b) + \n(a,b)(1,1) \in V_{n+1}.
    \]
    
    First put $a=(0,\dots,0,(\sqrt{2})^{n-2})$ and $b = z$. The inclusion $a \in M_n$ can be proved inductively by using Proposition~\ref{proposition:new-elements}. Then $\t(a) = 0$ and $\n(a,b) = (\sqrt{2})^{n-2}$, whence $(\sqrt{2})^{n-2}(1,1) \in V_{n+1}$.

    It follows that $\t(a)(b,b) \in V_{n+1}$ for all $a,b \in M_n$. Putting $a = z$, we obtain that $(b,b) \in M_{n+1}$ for an arbitrary $b \in M_n$.
\end{proof}

\section{Products of elements of ordinary Barnes--Wall lattices} \label{section:integer}

In this section we consider the embedding of the lattices~$L_n$ and~$L'_n$ into the Cayley--Dickson algebra~$\CD_n$ and study the sets of pairwise products of their elements. In particular, we are interested in whether~$L_n$ and~$L'_n$ are closed under multiplication.

\begin{notation}
    Let $X,Y \subset \CD_n$ be some lattices (over~$\Z$). Then $X \cdot Y$ is the set of integer linear combinations of elements of the form~$xy$, $x \in X$, $y \in Y$. In other words, $X \cdot Y \subset \CD_n$ is the lattice generated by all products whose first factor belongs to~$X$ and second factor belongs to~$Y$. 
\end{notation}

\begin{theorem} \label{theorem:classic-closed}
    For $n \geq 2$ the following inclusions hold:
    \[
    L_n \cdot L_n \subseteq L_n, \qquad L_n \cdot L'_n \subseteq L'_n, \qquad L'_n \cdot L_n \subseteq L'_n, \qquad L'_n \cdot L'_n \subseteq L_n/2.
    \]
\end{theorem}

\begin{proof}
    By definition, $M_n = L_n + \sqrt{2}L'_n$, and all elements of~$L_n$ and~$L'_n$ have rational coordinates. As follows from the definition of the Cayley--Dickson product, the product of elements with rational coordinates also has rational coordinates. According to Theorem~\ref{theorem:balanced-closed}, $M_n$ is closed under multiplication, hence we obtain
    \[
    M_n \supseteq M_n \cdot M_n = (L_n \cdot L_n + 2 L'_n \cdot L'_n) + \sqrt{2} (L_n \cdot L'_n + L'_n \cdot L_n),
    \]
    and thus 
    \[
        L_n \supseteq L_n \cdot L_n + 2 L'_n \cdot L'_n,\qquad
        L'_n \supseteq L_n \cdot L'_n + L'_n \cdot L_n.
    \]
    These inclusions immediately imply the required statement.
\end{proof}

\begin{remark} \label{remark:dimension-2}
    For $n = 2$ all inclusions in Theorem~\ref{theorem:classic-closed} turn into equalities. Indeed, in all cases except the last one this follows from the fact that $1 \in L_2$ is the first row of $H_2$. The equality $L'_2 \cdot L'_2 = L_2/2$ can be verified directly by dividing by $\sqrt{2}$ those rows and columns of the multiplication table~\ref{table:M2-multiplication} which correspond to the elements $(1+i)/\sqrt{2}$ and $(1+j)/\sqrt{2}$.
\end{remark}

\begin{proposition} \label{proposition:first-row}
    The first rows of the matrices $H_n$ and $H'_n$ are equal to $(2^k,0,\dots,0)$ and $(2^{k'},0,\dots,0)$, respectively, where $k = \lceil n/2 \rceil - 1$ and $k' = \lfloor n/2 \rfloor - 1$. Moreover, if an integer is divisible by $(\sqrt{2})^{n-2}$, then it is also divisible by $2^k$.
\end{proposition}

\begin{proof}
    By Proposition~\ref{proposition:main-properties}(4), the first row of the matrix $G_n$ equals $y = ((\sqrt{2})^{n-2}, 0, \dots, 0)$. If $n$ is even, then $k = k' = n/2 - 1$, the row $y$ has type~I, and the first rows of $H_n$ and~$H'_n$ coincide with it, see Proposition~\ref{proposition:common-types}. If $n$ is odd, then $y$ has type~II, so the first row of $H_n$ equals $((\sqrt{2})^{n-1}, 0, \dots, 0) = (2^k,0,\dots,0)$, while the first row of $H'_n$ equals $((\sqrt{2})^{n-3}, 0, \dots, 0) = (2^{k'},0,\dots,0)$.

    Note that $2^k$ is the smallest rational power of $\sqrt{2}$ which is divisible by $(\sqrt{2})^{n-2}$, and $2^{k'}$ is the largest rational power of $\sqrt{2}$ which divides $(\sqrt{2})^{n-2}$. This immediately implies the second part of the statement.
\end{proof}

\begin{proposition} \label{proposition:trace-classic}
    For all $x,y \in L_n$, $x',y' \in L'_n$, $n \geq 1$, we have 
    \[
    \t(x), \t(x'), \n(x,y) \in \Z,
    \]
    and also
    \[
        \n(x',y) \in \begin{cases}
            \Z/2, & n = 1,\\
            \Z, & n \geq 2;\\
        \end{cases} \qquad\quad
        \n(x',y') \in \begin{cases}
            \Z/2, & 1 \leq n \leq 2,\\
            \Z, & n \geq 3.\\
        \end{cases}
    \]
\end{proposition}

\begin{proof}
    Recall that $L_n, \sqrt{2} L'_n \subset M_n$. By Proposition~\ref{proposition:trace}, for all $z, w \in M_n$ we have $\t(z) \in \Z[\sqrt{2}]$, and~$\n(z,w)$ is divisible by $(\sqrt{2})^{n-2}$. Consequently, $\t(x) \in \Z[\sqrt{2}]$, $\t(x') \in \Z[\sqrt{2}]/\sqrt{2}$, and also~$\n(x,y)$, $\n(x',y)$, and~$\n(x',y')$ are divisible by $(\sqrt{2})^{n-2}$, $(\sqrt{2})^{n-3}$, and $(\sqrt{2})^{n-4}$, respectively. The elements with rational coordinates have rational values of the trace and the norm, which immediately implies the required statement.
\end{proof}

\begin{notation}
    In the proofs of Theorems~\ref{theorem:classic-generated} and~\ref{theorem:second-generated} we will occasionally use the notation $(a,b)_{n} = (a,b) \in \CD_{n}$, where $a,b \in \CD_{n-1}$. The index $n$ indicates explicitly the number of the current step of the Cayley--Dickson doubling process.
\end{notation}

The first half of the proof of the following theorem is similar to the proof of Theorem~\ref{theorem:balanced-generated}, with minor changes.

\begin{theorem} \label{theorem:classic-generated}
    If $n \geq 2$, then all pairwise products of elements of~$L_n$ generate~$L_n$ (over~$\Z$), that is, $L_n \cdot L_n = L_n$.
\end{theorem}

\begin{proof}
    The case $n = 2$ follows from Remark~\ref{remark:dimension-2}. For convenience of notation, we further prove the statement for~$L_{n+1}$, $n \geq 2$. Denote $V_{n+1} = L_{n+1} \cdot L_{n+1}$. It suffices to show that all rows of the generator matrix~$H_{n+1}$ belong to~$V_{n+1}$.

    Let $a, b \in L_n$ be arbitrary elements. By Corollaries~\ref{corollary:new-elements} and~\ref{corollary:new-elements-Ln}, we have $(a,a),(b,b),$ $(0,2a),(2b,0) \in L_{n+1}$. Consider the equality
    \[
        V_{n+1} \ni (b,b)(0,2a) + (a,a)(2b,0) = \t(b)(0,2a) - \t(a)(2b,0) + (4ab,0).
    \]
    Since $(4ab,0) = (2a,0)(2b,0) \in V_{n+1}$, we have 
    \begin{equation} \label{equation:generated-sum-classic}
        \t(b)(0,2a) - \t(a)(2b,0) \in V_{n+1}.
    \end{equation}
    
    Let $z = (1/2, \dots, 1/2)$ be the last row of $H_n$. Then $z$ has type I, $\t(z) = 1$, and
    \[
        V_{n+1} \ni (2z,0)(z,z) - (z,z)(0,2z) = \t(z)(2z,0) = (2z,0).
    \]
    Putting $a = b = z$ in Eq.~\eqref{equation:generated-sum-classic}, we obtain that $(0,2z) \in V_{n+1}$.

    Next, if $b = z$ and $a \in L_n$ is an arbitrary element, then $\t(b) = 1$, and Eq.~\eqref{equation:generated-sum-classic} implies $(0,2a) \in V_{n+1}$. Similarly, if $a = z$ and $b \in L_n$ is an arbitrary element, then $(2b,0) \in V_{n+1}$.

    For arbitrary $a,b \in L_n$ consider the equality
    \[
        V_{n+1} \ni (a,a)(b,b) = \t(a)(b,b) + \n(a,b)(1,1) - \t(a)(0,2b) + (0,2ba) - (2\n(a,b),0).
    \]
    By Proposition~\ref{proposition:first-row}, the first row of~$H_n$ equals $(2^k,0,\dots,0)$, where $k = \lceil n/2 \rceil - 1$. Besides, Theorem~\ref{theorem:classic-closed} implies $ba \in L_n$. Then, by what has been proved above, we have
    \[
    (0,2b), (0,2ba), (2 \cdot 2^k,0) \in V_{n+1}.
    \]
    By Propositions~\ref{proposition:trace}(2) and~\ref{proposition:trace-classic}, the value~$\n(a,b)$ is an integer which is divisible by $(\sqrt{2})^{n-2}$ in $\Z[\sqrt{2}]$. Applying Proposition~\ref{proposition:first-row}, we obtain that~$\n(a,b)$ is divisible by~$2^k$, hence $(2\n(a,b),0) \in V_{n+1}$. Therefore,
    \[
    \t(a)(b,b) + \n(a,b)(1,1) \in V_{n+1}.
    \]
    First put $a=(0,\dots,0,2^k)$ and $b = z$. Then $\t(a) = 0$ and $\n(a,b) = 2^k$, hence $2^k(1,1) \in V_{n+1}$. It follows that $\t(a)(b,b) \in V_{n+1}$ for all $a,b \in L_n$. Putting $a = z$, we obtain that $(b,b) \in L_{n+1}$ for an arbitrary $b \in L_n$.

    Thus, by Corollary~\ref{corollary:new-elements-types}, it remains to show that~$V_{n+1}$ contains the elements of the form~$(a,0)$, where~$a$ is a row of~$H_n$ of type II. For an arbitrary $b \in L_n$ consider the equality
    \begin{align*}
        V_{n+1} \ni (a,0)(b,b) &= (ab,ba) = (\t(b)a + \t(a)b - \n(a,b) - ba,ba) \\
        &= \t(b)(a,0) + \t(a)(b,0) - (\n(a,b),0) + (ba,ba) - (2ba,0).
    \end{align*}
    By what has been proved above, $(2b,0), (ba,ba), (2ba,0) \in V_{n+1}$. The row $a$ has type II, hence it is not the last one, and therefore $\t(a)$ is divisible by $2$. Consequently,
    \begin{equation} \label{equation:generating-type-II}
        \t(b)(a,0) - (\n(a,b),0) \in V_{n+1}.
    \end{equation}
    
    Since the row $a$ has type II, it is $\sqrt{2}$ times greater than the corresponding row $a'$ of the matrix~$G_n$. By Proposition~\ref{proposition:rows-norm}, $\n(a',z)$ is divisible by $(\sqrt{2})^{n-2}$, and the equality is attained if and only if $a'$ is the first row of $G_n$, otherwise $\n(a',z)$ exceeds $(\sqrt{2})^{n-2}$ by a factor of at least $\sqrt{2}$.
    
    In the second case $\n(a,z)$ is an integer which exceeds $(\sqrt{2})^{n-2}$ by a factor of at least $2$, and therefore $\n(a,z)$ is divisible by $2 \cdot 2^k$. Then, as was proved above, $(\n(a,z),0) \in V_{n+1}$. Putting $b = z$ in Eq.~\eqref{equation:generating-type-II}, we obtain that $(a,0) \in V_{n+1}$.

    It remains to show that if the first row $y$ of the matrix $H_n$ has type II, then $(y,0) \in V_{n+1}$. In this case $y = (2^k,0,\dots,0) = (2\widetilde{y},0)_n$ and $z = (\widetilde{z},\widetilde{z})_n$, where $\widetilde{y}$ and $\widetilde{z}$ are the first and the last rows of the matrix $H_{n-1}$, respectively, and both of them have type~I. Put $a = (2\widetilde{z},0)_n$ and $b = (\widetilde{y},\widetilde{y})_n$. Then $\t(b) = \t(\widetilde{y})$ is divisible by $2$, since $n \geq 2$ is an odd integer, and 
    \[
    \n(a,b) = \n((2\widetilde{z},0)_n,(\widetilde{y},\widetilde{y})_n) = 2\n(\widetilde{z},\widetilde{y}) = \n((2\widetilde{y},0)_n,(\widetilde{z},\widetilde{z})_n) = \n(y,z) = 2^k.
    \]
    The row~$a$ has type~II, and $(2a,0) \in V_{n+1}$ by what has been proved above. Hence Eq.~\eqref{equation:generating-type-II} yields $(y,0) = (2^k,0) \in V_{n+1}$.
\end{proof}

\begin{theorem} \label{theorem:second-generated}
    If $n \geq 2$, then $L_n \cdot L'_n = L'_n$.
\end{theorem}

\begin{proof}
    The case $n = 2$ follows from Remark~\ref{remark:dimension-2}. Denote $V_{n+1} = L_{n+1} \cdot L'_{n+1}$. It suffices to show that for $n \geq 2$ all rows of the generator matrix~$H'_{n+1}$ belong to~$V_{n+1}$. By Proposition~\ref{proposition:connection-two-lattices}, they have the form~$(a,0)$ and~$(b',b')$, where~$a$ and~$b'$ are arbitrary rows of the matrices~$H_n$ and~$H'_n$, respectively.

    The chain of inclusions~\eqref{equation:inclusion} gives $L_{n+1} \subset L'_{n+1}$, so Theorem~\ref{theorem:classic-generated} implies 
    \[
    L_{n+1} = L_{n+1} \cdot L_{n+1} \subset L_{n+1} \cdot L'_{n+1} = V_{n+1}.
    \]
    Let $z = (1/2, \dots, 1/2)$ be the last row of~$H_n$ and~$H'_n$. Then $\t(z) = 1$, and Corollary~\ref{corollary:new-elements} yields $(z,z) \in L_{n+1}$ and $(0,z) \in L'_{n+1}$, whence
    \[
        V_{n+1} \ni (z,z)(0,z) = (-\bar{z}z,z^2) = (-\t(z)z + z^2,z^2) = -(z,0) + (z^2,z^2).
    \]
    By Theorem~\ref{theorem:classic-closed}, $z^2 \in L_n$, hence $(z^2,z^2) \in L_{n+1} \subset V_{n+1}$. Consequently, $(z,0) \in V_{n+1}$.
    
    Let $a, b \in L_n$ be arbitrary elements. In what follows we use Corollary~\ref{corollary:new-elements} without further explanations. Consider the equality
    \[
        V_{n+1} \ni (b,b)(0,a) + (a,a)(b,0) = \t(b)(0,a) - \t(a)(b,0) + (2ab,0).
    \]
    Since $(2ab,0) = (2a,0)(b,0) \in V_{n+1}$, we have 
    \[
        \t(b)(0,a) - \t(a)(b,0) \in V_{n+1}.
    \]
    If $b = z$ and $a \in L_n$ is an arbitrary element, then $\t(b) = 1$ and $(0,a) \in V_{n+1}$. Similarly, if $a = z$ and $b \in L_n$ is an arbitrary element, then $(b,0) \in V_{n+1}$.

    Further, for arbitrary $a \in L_n$ and $b' \in L'_n$ consider the equality
    \[
        V_{n+1} \ni (a,a)(b',b') = \t(a)(b',b') + \n(a,b')(1,1) - (2\n(a,b'),0) - \t(a)(0,2b') + (0,2b'a).
    \]
    By Theorem~\ref{theorem:classic-closed}, $b'a \in L'_n$, hence $(0,2b'),(0,2b'a) \in L_{n+1} \subset V_{n+1}$.
    Thus,
    \begin{equation} \label{equation:generating-second}
    \t(a)(b',b') + \n(a,b')(1,1) - (2\n(a,b'),0) \in V_{n+1}.
    \end{equation}

    It remains to show that $(b',b') \in V_{n+1}$, where $b'$ is an arbitrary row of the matrix~$H'_n$. If $b'$ has type~I, then, by Proposition~\ref{proposition:common-types}, it coincides with the corresponding row of $H_n$, and then the required fact follows from $(b',b') \in L_{n+1} \subset V_{n+1}$. Otherwise the row $b'$ has type~II, and it is $\sqrt{2}$ times smaller than the corresponding row $b$ of the matrix $G_n$. By Proposition~\ref{proposition:rows-norm}, $\n(z,b)$ is divisible by $(\sqrt{2})^{n-2}$, and the equality is attained if and only if $b$ is the first row of $H_n$, otherwise $\n(z,b)$ exceeds $(\sqrt{2})^{n-2}$ by a factor of at least $\sqrt{2}$.
    
    In the second case, by Proposition~\ref{proposition:trace-classic}, the value $\n(z,b')$ is an integer divisible by $(\sqrt{2})^{n-2}$, and therefore $\n(z,b')$ is divisible by $2^k$, where $k = \lceil n/2 \rceil - 1$, see Proposition~\ref{proposition:first-row}. The first row of~$H_n$ equals $(2^k,0,\dots,0)$, hence $(\n(z,b'), 0, \dots, 0) \in L_n$, and consequently
    \[
    \n(z,b')(1,1), (2\n(z,b'),0) \in L_{n+1} \subset V_{n+1}.
    \]
    Putting $a = z$ in Eq.~\eqref{equation:generating-second}, we obtain that $(b',b') \in V_{n+1}$.

    Finally, we show that if the first row~$y$ of the matrix~$H'_n$ has type~II, then $(y,y) \in V_{n+1}$. In this case it follows from Propositions~\ref{proposition:common-types}, \ref{proposition:connection-two-lattices}, and~\ref{proposition:first-row} that $y = (2^{k-1},0,\dots,0) = (\widetilde{y},0)_n$ and $z = (\widetilde{z},\widetilde{z})_n$, where $\widetilde{y}$ and $\widetilde{z}$ are the first and the last rows of the matrix~$H_{n-1}$, respectively.
    Put $a = (\widetilde{y},\widetilde{y})_n \in L_n$ and $b' = (\widetilde{z},0)_n \in L'_n$. Then $\t(a) = \t(\widetilde{y})$ is divisible by~$2$, since $n \geq 2$ is an odd integer, and 
    \[
    \n(a,b') = \n((\widetilde{y},\widetilde{y})_n,(\widetilde{z},0)_n) = \n(\widetilde{z},\widetilde{y}) = \n((\widetilde{z},\widetilde{z})_n,(\widetilde{y},0)_n) = \n(z,y) = 2^{k-1}.
    \]
    It follows from $2b' \in 2L'_{n} \subset L_n$ that $2(b',b') \in L_{n+1} \subset V_{n+1}$. Besides, by what has been proved above, $(2^k,0,\dots,0) \in L_n$ implies $(2\n(a,b'),0) = (2^k,0) \in V_{n+1}$. Hence Eq.~\eqref{equation:generating-second} yields $(y,y) = 2^{k-1}(1,1) \in V_{n+1}$.
\end{proof}

    \begin{remark} \label{remark:strict-inclusion}
        For $n \geq 3$ the inclusions $L'_n \cdot L_n \subset L'_n$ and $L'_n \cdot L'_n \subset L_n/2$ from Theorem~\ref{theorem:classic-closed} can be strict. Denote the generator matrices of the lattices $L'_n \cdot L_n$ and $L'_n \cdot L'_n$ by~$C_n$ and~$C'_n$, respectively. Then for $n = 3$ the lower triangular Hermite normal forms of the generator matrices are as follows:
    \begin{center}
    \resizebox{0.97\textwidth}{!}{$
        C_3 \to \widetilde{C}_3 = \left(\begin{array}{cccccccc}
        2 & 0 & 0 & 0 & 0 & 0 & 0 & 0 \\
        1 & 1 & 0 & 0 & 0 & 0 & 0 & 0 \\
        1 & 0 & 1 & 0 & 0 & 0 & 0 & 0 \\
        1 & 0 & 0 & 1 & 0 & 0 & 0 & 0 \\
        1 & 0 & 0 & 0 & 1 & 0 & 0 & 0\\
        1 & 0 & 0 & 0 & 0 & 1 & 0 & 0\\
        1 & 0 & 0 & 0 & 0 & 0 & 1 & 0\\
        \frac{1}{2} & \frac{1}{2} & \frac{1}{2} & \frac{1}{2} & \frac{1}{2} & \frac{1}{2} & \frac{1}{2} & \frac{1}{2} \\
        \end{array}\right), \quad
        C'_3 \to \widetilde{C}'_3 = \left(\begin{array}{cccccccc}
         1 & 0 & 0 & 0 & 0 & 0 & 0 & 0 \\
         \frac{1}{2} & \frac{1}{2} & 0 & 0 & 0 & 0 & 0 & 0 \\
         \frac{1}{2} & 0 & \frac{1}{2} & 0 & 0 & 0 & 0 & 0 \\
         \frac{1}{2} & 0 & 0 & \frac{1}{2} & 0 & 0 & 0 & 0 \\
         \frac{1}{2} & 0 & 0 & 0 & \frac{1}{2} & 0 & 0 & 0 \\
         \frac{1}{2} & 0 & 0 & 0 & 0 & \frac{1}{2} & 0 & 0 \\
         \frac{1}{2} & 0 & 0 & 0 & 0 & 0 & \frac{1}{2} & 0 \\
        \frac{1}{2} & 0 & 0 & 0 & 0 & 0 & 0 & \frac{1}{2} \\
        \end{array}\right).
    $}
    \end{center}
    Recall that the matrices $\widetilde{H}_3$ and $\widetilde{H}'_3$ were computed in Example~\ref{example:low-dimension}.
    We also write down the Hermite normal forms for $n = 4$:

    \begin{center}
    \resizebox{0.75\textwidth}{!}{$
    \begin{aligned}
        H_4 &\to \widetilde{H}_4 = \left(
        \begin{array}{cccccccccccccccc}
         2 & 0 & 0 & 0 & 0 & 0 & 0 & 0 & 0 & 0 & 0 & 0 & 0 & 0 & 0 & 0 \\
         0 & 2 & 0 & 0 & 0 & 0 & 0 & 0 & 0 & 0 & 0 & 0 & 0 & 0 & 0 & 0 \\
         0 & 0 & 2 & 0 & 0 & 0 & 0 & 0 & 0 & 0 & 0 & 0 & 0 & 0 & 0 & 0 \\
         1 & 1 & 1 & 1 & 0 & 0 & 0 & 0 & 0 & 0 & 0 & 0 & 0 & 0 & 0 & 0 \\
         0 & 0 & 0 & 0 & 2 & 0 & 0 & 0 & 0 & 0 & 0 & 0 & 0 & 0 & 0 & 0 \\
         1 & 1 & 0 & 0 & 1 & 1 & 0 & 0 & 0 & 0 & 0 & 0 & 0 & 0 & 0 & 0 \\
         1 & 0 & 1 & 0 & 1 & 0 & 1 & 0 & 0 & 0 & 0 & 0 & 0 & 0 & 0 & 0 \\
         0 & 1 & 1 & 0 & 1 & 0 & 0 & 1 & 0 & 0 & 0 & 0 & 0 & 0 & 0 & 0 \\
         0 & 0 & 0 & 0 & 0 & 0 & 0 & 0 & 2 & 0 & 0 & 0 & 0 & 0 & 0 & 0 \\
         1 & 1 & 0 & 0 & 0 & 0 & 0 & 0 & 1 & 1 & 0 & 0 & 0 & 0 & 0 & 0 \\
         1 & 0 & 1 & 0 & 0 & 0 & 0 & 0 & 1 & 0 & 1 & 0 & 0 & 0 & 0 & 0 \\
         0 & 1 & 1 & 0 & 0 & 0 & 0 & 0 & 1 & 0 & 0 & 1 & 0 & 0 & 0 & 0 \\
         1 & 0 & 0 & 0 & 1 & 0 & 0 & 0 & 1 & 0 & 0 & 0 & 1 & 0 & 0 & 0 \\
         0 & 1 & 0 & 0 & 1 & 0 & 0 & 0 & 1 & 0 & 0 & 0 & 0 & 1 & 0 & 0 \\
         0 & 0 & 1 & 0 & 1 & 0 & 0 & 0 & 1 & 0 & 0 & 0 & 0 & 0 & 1 & 0 \\
         \frac{1}{2} & \frac{1}{2} & \frac{1}{2} & \frac{1}{2} & \frac{1}{2} & \frac{1}{2} & \frac{1}{2} & \frac{1}{2} & \frac{1}{2} & \frac{1}{2} & \frac{1}{2} & \frac{1}{2} & \frac{1}{2} & \frac{1}{2} &
           \frac{1}{2} & \frac{1}{2} \\
        \end{array}
        \right),\\
        H'_4 &\to \widetilde{H}'_4 = \left(
        \begin{array}{cccccccccccccccc}
         2 & 0 & 0 & 0 & 0 & 0 & 0 & 0 & 0 & 0 & 0 & 0 & 0 & 0 & 0 & 0 \\
         1 & 1 & 0 & 0 & 0 & 0 & 0 & 0 & 0 & 0 & 0 & 0 & 0 & 0 & 0 & 0 \\
         1 & 0 & 1 & 0 & 0 & 0 & 0 & 0 & 0 & 0 & 0 & 0 & 0 & 0 & 0 & 0 \\
         1 & 0 & 0 & 1 & 0 & 0 & 0 & 0 & 0 & 0 & 0 & 0 & 0 & 0 & 0 & 0 \\
         1 & 0 & 0 & 0 & 1 & 0 & 0 & 0 & 0 & 0 & 0 & 0 & 0 & 0 & 0 & 0 \\
         1 & 0 & 0 & 0 & 0 & 1 & 0 & 0 & 0 & 0 & 0 & 0 & 0 & 0 & 0 & 0 \\
         1 & 0 & 0 & 0 & 0 & 0 & 1 & 0 & 0 & 0 & 0 & 0 & 0 & 0 & 0 & 0 \\
         \frac{1}{2} & \frac{1}{2} & \frac{1}{2} & \frac{1}{2} & \frac{1}{2} & \frac{1}{2} & \frac{1}{2} & \frac{1}{2} & 0 & 0 & 0 & 0 & 0 & 0 & 0 & 0 \\
         1 & 0 & 0 & 0 & 0 & 0 & 0 & 0 & 1 & 0 & 0 & 0 & 0 & 0 & 0 & 0 \\
         1 & 0 & 0 & 0 & 0 & 0 & 0 & 0 & 0 & 1 & 0 & 0 & 0 & 0 & 0 & 0 \\
         1 & 0 & 0 & 0 & 0 & 0 & 0 & 0 & 0 & 0 & 1 & 0 & 0 & 0 & 0 & 0 \\
         \frac{1}{2} & \frac{1}{2} & \frac{1}{2} & \frac{1}{2} & 0 & 0 & 0 & 0 & \frac{1}{2} & \frac{1}{2} & \frac{1}{2} & \frac{1}{2} & 0 & 0 & 0 & 0 \\
         1 & 0 & 0 & 0 & 0 & 0 & 0 & 0 & 0 & 0 & 0 & 0 & 1 & 0 & 0 & 0 \\
         \frac{1}{2} & \frac{1}{2} & 0 & 0 & \frac{1}{2} & \frac{1}{2} & 0 & 0 & \frac{1}{2} & \frac{1}{2} & 0 & 0 & \frac{1}{2} & \frac{1}{2} & 0 & 0 \\
         \frac{1}{2} & 0 & \frac{1}{2} & 0 & \frac{1}{2} & 0 & \frac{1}{2} & 0 & \frac{1}{2} & 0 & \frac{1}{2} & 0 & \frac{1}{2} & 0 & \frac{1}{2} & 0 \\
         \frac{1}{2} & 0 & 0 & \frac{1}{2} & 0 & \frac{1}{2} & \frac{1}{2} & 0 & 0 & \frac{1}{2} & \frac{1}{2} & 0 & \frac{1}{2} & 0 & 0 & \frac{1}{2} \\
        \end{array}
        \right),\\
    \end{aligned}
    $}
    \end{center}

    \begin{center}
    \resizebox{0.75\textwidth}{!}{$
    \begin{aligned}
        C_4 &\to \widetilde{C}_4 = \left(
        \begin{array}{cccccccccccccccc}
         2 & 0 & 0 & 0 & 0 & 0 & 0 & 0 & 0 & 0 & 0 & 0 & 0 & 0 & 0 & 0 \\
         1 & 1 & 0 & 0 & 0 & 0 & 0 & 0 & 0 & 0 & 0 & 0 & 0 & 0 & 0 & 0 \\
         1 & 0 & 1 & 0 & 0 & 0 & 0 & 0 & 0 & 0 & 0 & 0 & 0 & 0 & 0 & 0 \\
         1 & 0 & 0 & 1 & 0 & 0 & 0 & 0 & 0 & 0 & 0 & 0 & 0 & 0 & 0 & 0 \\
         1 & 0 & 0 & 0 & 1 & 0 & 0 & 0 & 0 & 0 & 0 & 0 & 0 & 0 & 0 & 0 \\
         1 & 0 & 0 & 0 & 0 & 1 & 0 & 0 & 0 & 0 & 0 & 0 & 0 & 0 & 0 & 0 \\
         1 & 0 & 0 & 0 & 0 & 0 & 1 & 0 & 0 & 0 & 0 & 0 & 0 & 0 & 0 & 0 \\
         1 & 0 & 0 & 0 & 0 & 0 & 0 & 1 & 0 & 0 & 0 & 0 & 0 & 0 & 0 & 0 \\
         1 & 0 & 0 & 0 & 0 & 0 & 0 & 0 & 1 & 0 & 0 & 0 & 0 & 0 & 0 & 0 \\
         1 & 0 & 0 & 0 & 0 & 0 & 0 & 0 & 0 & 1 & 0 & 0 & 0 & 0 & 0 & 0 \\
         1 & 0 & 0 & 0 & 0 & 0 & 0 & 0 & 0 & 0 & 1 & 0 & 0 & 0 & 0 & 0 \\
         1 & 0 & 0 & 0 & 0 & 0 & 0 & 0 & 0 & 0 & 0 & 1 & 0 & 0 & 0 & 0 \\
         1 & 0 & 0 & 0 & 0 & 0 & 0 & 0 & 0 & 0 & 0 & 0 & 1 & 0 & 0 & 0 \\
         1 & 0 & 0 & 0 & 0 & 0 & 0 & 0 & 0 & 0 & 0 & 0 & 0 & 1 & 0 & 0 \\
         1 & 0 & 0 & 0 & 0 & 0 & 0 & 0 & 0 & 0 & 0 & 0 & 0 & 0 & 1 & 0 \\
         \frac{1}{2} & \frac{1}{2} & \frac{1}{2} & \frac{1}{2} & \frac{1}{2} & \frac{1}{2} & \frac{1}{2} & \frac{1}{2} & \frac{1}{2} & \frac{1}{2} & \frac{1}{2} & \frac{1}{2} & \frac{1}{2} & \frac{1}{2} &
           \frac{1}{2} & \frac{1}{2} \\
        \end{array}
        \right),\\
        C'_4 &\to \widetilde{C}'_4 = \left(
        \begin{array}{cccccccccccccccc}
         1 & 0 & 0 & 0 & 0 & 0 & 0 & 0 & 0 & 0 & 0 & 0 & 0 & 0 & 0 & 0 \\
         0 & 1 & 0 & 0 & 0 & 0 & 0 & 0 & 0 & 0 & 0 & 0 & 0 & 0 & 0 & 0 \\
         0 & 0 & 1 & 0 & 0 & 0 & 0 & 0 & 0 & 0 & 0 & 0 & 0 & 0 & 0 & 0 \\
         \frac{1}{2} & \frac{1}{2} & \frac{1}{2} & \frac{1}{2} & 0 & 0 & 0 & 0 & 0 & 0 & 0 & 0 & 0 & 0 & 0 & 0 \\
         0 & 0 & 0 & 0 & 1 & 0 & 0 & 0 & 0 & 0 & 0 & 0 & 0 & 0 & 0 & 0 \\
         \frac{1}{2} & \frac{1}{2} & 0 & 0 & \frac{1}{2} & \frac{1}{2} & 0 & 0 & 0 & 0 & 0 & 0 & 0 & 0 & 0 & 0 \\
         \frac{1}{2} & 0 & \frac{1}{2} & 0 & \frac{1}{2} & 0 & \frac{1}{2} & 0 & 0 & 0 & 0 & 0 & 0 & 0 & 0 & 0 \\
         0 & \frac{1}{2} & \frac{1}{2} & 0 & \frac{1}{2} & 0 & 0 & \frac{1}{2} & 0 & 0 & 0 & 0 & 0 & 0 & 0 & 0 \\
         0 & 0 & 0 & 0 & 0 & 0 & 0 & 0 & 1 & 0 & 0 & 0 & 0 & 0 & 0 & 0 \\
         0 & 0 & 0 & 0 & 0 & 0 & 0 & 0 & 0 & 1 & 0 & 0 & 0 & 0 & 0 & 0 \\
         0 & 0 & 0 & 0 & 0 & 0 & 0 & 0 & 0 & 0 & 1 & 0 & 0 & 0 & 0 & 0 \\
         0 & 0 & 0 & 0 & 0 & 0 & 0 & 0 & \frac{1}{2} & \frac{1}{2} & \frac{1}{2} & \frac{1}{2} & 0 & 0 & 0 & 0 \\
         0 & 0 & 0 & 0 & 0 & 0 & 0 & 0 & 0 & 0 & 0 & 0 & 1 & 0 & 0 & 0 \\
         0 & 0 & 0 & 0 & 0 & 0 & 0 & 0 & \frac{1}{2} & \frac{1}{2} & 0 & 0 & \frac{1}{2} & \frac{1}{2} & 0 & 0 \\
         0 & 0 & 0 & 0 & 0 & 0 & 0 & 0 & \frac{1}{2} & 0 & \frac{1}{2} & 0 & \frac{1}{2} & 0 & \frac{1}{2} & 0 \\
         0 & 0 & 0 & 0 & 0 & 0 & 0 & 0 & 0 & \frac{1}{2} & \frac{1}{2} & 0 & \frac{1}{2} & 0 & 0 & \frac{1}{2} \\
        \end{array}
        \right).
    \end{aligned}
    $}
    \end{center}
    
    In both cases $\widetilde{C}_n \neq \widetilde{H}'_n$ and $\widetilde{C}'_n \neq \widetilde{H}_n/2$, hence $L'_n \cdot L_n \subsetneq L'_n$ and $L'_n \cdot L'_n \subsetneq L_n/2$. 
    \end{remark}

    \begin{corollary}
        It follows from Remarks~\ref{remark:not-closed-dimension-one}, \ref{remark:dimension-2}, and~\ref{remark:strict-inclusion} that for $1 \leq n \leq 4$ the lattice~$L'_n$ is not closed under multiplication.
    \end{corollary}

    The following proposition shows how the lattices $L_n$ and~$L'_n$ are related to each other in the Cayley--Dickson algebra $\CD_n$.
    
    \begin{proposition} \label{proposition:connection}
        For all $n \in \N$ we have
        \[
        L_n = L'_n (1 + e_{2^{n-1}}) \qquad \text{and} \qquad L'_n = L_n (1 + e_{2^{n-1}})/2.
        \]
    \end{proposition}
    
    \begin{proof}
        By Proposition~\ref{proposition:connection-two-lattices}, the rows of~$H_{n}$ have the form~$(x,x)$ and~$(2x',0)$, while the rows of~$H'_{n}$ have the form~$(x',x')$ and~$(x,0)$, where~$x$ is an arbitrary row of~$H_{n-1}$, and~$x'$ is an arbitrary row of~$H'_{n-1}$. Moreover, $1 + e_{2^{n-1}} = (1,1) \in \CD_{n}$. Then
        \begin{align*}
            (x,x)(1,1) &= (0,2x), & (2x',0)(1,1) &= (2x',2x'),\\
            (x',x')(1,1) &= (0,2x'), & (x,0)(1,1) &= (x,x).        
        \end{align*}
        Thus the rows of the generator matrix of the lattice $L'_n (1 + e_{2^{n-1}})$ have the form $(x,x)$ and $(0,2x')$, and the rows of the generator matrix of $L_n (1 + e_{2^{n-1}})/2$ have the form $(x',x')$ and $(0,x)$. Taking into account the chain of inclusions~\eqref{equation:inclusion} and Corollary~\ref{corollary:new-elements}, we obtain the required equalities of lattices.
    \end{proof}
    
    \begin{remark} 
        Instead of the standard basis for the octonion algebra~$\mathbb{O}$, Wilson~\cite{Wilson} used the basis which consists of the elements $1,i_0,\dots,i_6$ such that $i_0 i_{1} = -i_{1}i_0 = i_{3}$, and this equality is preserved under the permutations of indices $t \mapsto t + 1 \pmod 7$ and $t \mapsto 2t \pmod 7$. He considered the lattice~$B$ generated by the elements of this basis and by the octonions of the form $(\pm 1 \pm i_t \pm i_{t+1} \pm i_{t+3})/2$, as well as by their ``complements'' $(\pm i_{t+2} \pm i_{t+4} \pm i_{t+5} \pm i_{t+6})/2$. The elements of this lattice are exactly the Kirmse integers, and it is well known that they are not closed under multiplication. Then Wilson put 
        \[
        A = (1+i_0)B(1+i_0)/2, \quad R = (1+i_0)B, \quad L = B(1+i_0),
        \]
        where~$A$ is one of the copies of the Coxeter--Dickson integral octonions closed under multiplication, $R$ corresponds to the lattice~$E_8$ written in the standard form via the even coordinate system, and~$L$ differs from~$R$ by the change of sign of one coordinate (for example, the last one). Using the Moufang identities, which hold in any alternative algebra, Wilson showed that $B \cdot L = L$ and $R \cdot B = R$.
    
        In terms of the standard basis, one can take the elements $1,e_6,e_1,e_4,e_7,e_5,e_3,-e_2$ as $1,i_0,\dots,i_6$. Then the lattice~$B$ coincides with~$L'_3$ and, due to the change of sign of the coordinate at~$i_6$, the lattice~$L$ coincides with~$L_3$, that is, $L_3 = L'_3 (1 + e_6)$. A similar relation between~$L_n$ and~$L'_n$ holds for arbitrary values of~$n$, see Proposition~\ref{proposition:connection}. Thus Wilson's first formula can be rewritten as $L'_3 \cdot L_3 = L_3$. It can be seen from Remark~\ref{remark:strict-inclusion} that $\widetilde{C}_3 = \widetilde{H}_3$, that is, the required equality indeed holds, but for $n = 4$ the inclusion $L'_4 \cdot L_4 \subset L_4$ is no longer true. This is explained by the fact that for $n \geq 4$ the algebras~$\CD_n$ are not alternative, hence the Moufang identities are not satisfied.
    \end{remark}

\medskip

The author is grateful to Academician of the Russian Academy of Sciences D.~O.~Orlov for drawing attention to this problem and to her scientific advisor Professor A.~E.~Guterman for his attention to this work and fruitful discussions.

\bibliographystyle{abbrvurl}
\bibliography{bibfile}

@article{Baez,
  author    = {John C. Baez},
  title     = {The Octonions},
  journal   = {Bull. Amer. Math. Soc.},
  volume    = {39},
  pages     = {145--205},
  year      = {2002},
  doi       = {10.1090/S0273-0979-01-00934-X}
}

@article{BarnesWall,
  author    = {Barnes, E. S. and Wall, G. E.},
  title     = {Some extreme forms defined in terms of abelian groups},
  journal   = {J.~Australian Math. Soc.},
  volume    = {1},
  pages     = {47--63},
  year      = {1959},
  doi       = {10.1017/S144678870000006X}
}

@article{Biss,
  author    = {Daniel K. Biss and Daniel Dugger and Daniel C. Isaksen},
  title     = {Large annihilators in {C}ayley--{D}ickson algebras},
  journal   = {Comm. Algebra},
  volume    = {36},
  number    = {2},
  pages     = {632--664},
  year      = {2008},
  doi       = {10.1080/00927870701724094}
}

@article{Biss2,
  author    = {Daniel K. Biss and Daniel Dugger and Daniel C. Isaksen},
  title     = {Large annihilators in {C}ayley--{D}ickson algebras~{II}},
  journal   = {Bol. Soc. Mat. Mex.},
  volume    = {13},
  number    = {2},
  pages     = {269--292},
  year      = {2007}
}

@article{Viazovska24,
  title     = {The sphere packing problem in dimension $24$},
  author    = {Cohn, Henry and Kumar, Abhinav and Miller, Stephen D. and Radchenko, Danylo and Viazovska, Maryna},
  journal   = {Annals of Mathematics},
  volume    = {185},
  number    = {3},
  pages     = {1017--1033},
  year      = {2017},
  publisher = {Princeton University and Institute for Advanced Study},
  doi       = {10.4007/annals.2017.185.3.8}
}

@book{ConwaySloane,
  title     = {Sphere Packings, Lattices and Groups},
  author    = {Conway, John H. and Sloane, Neil J. A.},
  edition   = {3rd},
  series    = {Grundlehren der mathematischen Wissenschaften},
  volume    = {290},
  year      = {1999},
  publisher = {Springer-Verlag},
  address   = {New York},
  doi       = {10.1007/978-1-4757-6568-7}
}

@article{ConwaySloane1982,
  author    = {Conway, J. H. and Sloane, N. J. A.},
  title     = {Twenty-three constructions for the {L}eech lattice},
  journal   = {Proc. R. Soc. Lond., Ser. A},
  volume    = {381},
  pages     = {275--283},
  year      = {1982},
  doi       = {10.1098/rspa.1982.0070}
}

@book{ConwaySmith,
  author    = {Conway, John H. and Smith, Derek A.},
  title     = {On Quaternions and Octonions: Their Geometry, Arithmetic, and Symmetry},
  year      = {2003},
  publisher = {A K Peters},
  address   = {Natick, Massachusetts},
  doi       = {10.1201/9781439864180}
}

@book{Dixon1994Book,
  title     = {Division Algebras: Octonions, Quaternions, Complex Numbers and the Algebraic Design of Physics},
  author    = {Dixon, Geoffrey M.},
  series    = {Mathematics and Its Applications},
  volume    = {290},
  year      = {1994},
  publisher = {Kluwer Academic Publishers},
  address   = {Dordrecht, Netherlands},
  doi       = {10.1007/978-1-4757-2315-1}
}

@book{Dixon2011Book,
  title     = {Division Algebras, Lattices, Physics, Windmill Tilting},
  author    = {Dixon, Geoffrey M.},
  year      = {2011},
  publisher = {CreateSpace Independent Publishing Platform},
  address   = {Charleston, SC}
}

@article{ForneyCosetII,
  author    = {Forney, Jr., G. David},
  title     = {Coset codes --- {P}art {II}. {B}inary lattices and related codes},
  journal   = {IEEE Transactions on Information Theory},
  year      = {1988},
  volume    = {34},
  number    = {5},
  pages     = {1152--1187},
  doi       = {10.1109/18.21246}
}

@misc{ForneyPrinciples,
  author       = {Forney, Jr., G. David},
  title        = {{Principles of Digital Communication II} {C}ourse Notes},
  year         = {2005},
  howpublished = {Massachusetts Institute of Technology},
  url          = {https://ocw.mit.edu/courses/6-451-principles-of-digital-communication-ii-spring-2005/resources/mit6_451s05_fulllecnotes/}
}

@article{Krausshar,
  author    = {Krau\ss{}har, Soören A.},
  title     = {Function Theories in {C}ayley--{D}ickson Algebras and Number Theory},
  journal   = {Milan J.~Math.},
  volume    = {89},
  pages     = {19--44},
  year      = {2021},
  doi       = {10.1007/s00032-021-00325-y}
}

@inproceedings{Micciancio,
  author    = {Micciancio, Daniele and Nicolosi, Antonio},
  booktitle = {2008 IEEE International Symposium on Information Theory},
  title     = {Efficient bounded distance decoders for {B}arnes--{W}all lattices},
  year      = {2008},
  pages     = {2484--2488},
  doi       = {10.1109/ISIT.2008.4595438}
}

@article{Moreno,
  author    = {Guillermo Moreno},
  title     = {The zero divisors of the {C}ayley--{D}ickson algebras over the real numbers},
  journal   = {Bol. Soc. Mat. Mex.},
  volume    = {4},
  number    = {1},
  pages     = {13--28},
  year      = {1998}
}

@incollection{Moreno_alternative,
  author    = {Guillermo Moreno},
  title     = {Alternative elements in the {C}ayley--{D}ickson algebras},
  booktitle = {Topics in Mathematical Physics, General Relativity and Cosmology in Honor of Jerzy Pleba\~{n}ski},
  pages     = {333--346},
  year      = {2006},
  publisher = {World Sci. Publ.},
  address   = {Hackensack, NJ},
  doi       = {10.1142/9789812772732_0028}
}

@misc{Moreno_constructing,
  author       = {Guillermo Moreno},
  title        = {Constructing zero divisors in the higher dimensional {C}ayley--{D}ickson algebras},
  archivePrefix= {arXiv},
  eprint       = {math/0512517},
  year         = {2005}
}

@article{Nebe1998,
  author    = {Nebe, Gabriele},
  title     = {Some cyclo-quaternionic lattices},
  journal   = {J.~Algebra},
  volume    = {199},
  number    = {2},
  pages     = {472--498},
  year      = {1998},
  doi       = {10.1006/jabr.1997.7163}
}

@incollection{Nebe1999,
  author    = {Nebe, Gabriele},
  title     = {Construction and investigation of lattices with matrix groups},
  booktitle = {Integral Quadratic Forms and Lattices (Seoul, 1998)},
  series    = {Contemporary Mathematics},
  volume    = {249},
  pages     = {205--220},
  publisher = {American Mathematical Society},
  address   = {Providence, RI},
  year      = {1999},
  doi       = {10.1090/conm/249/03759}
}

@incollection{Nebe2002,
  author    = {Nebe, G. and Rains, E. M. and Sloane, N. J. A.},
  editor    = {Blahut, Richard E. and Koetter, Ralf},
  title     = {A simple construction for the {B}arnes--{W}all lattices},
  booktitle = {Codes, graphs, and systems: A celebration of the life and career of {G}. {D}avid {F}orney, {J}r. on the occasion of his sixtieth birthday},
  year      = {2002},
  publisher = {Springer US},
  address   = {Boston, MA},
  pages     = {333--342},
  doi       = {10.1007/978-1-4615-0895-3_19}
}

@article{Quebbemann,
  title     = {Lattices with theta functions for $G(\sqrt{2})$ and linear codes},
  journal   = {J.~Algebra},
  volume    = {105},
  number    = {2},
  pages     = {443--450},
  year      = {1987},
  doi       = {10.1016/0021-8693(87)90208-0},
  author    = {H.-G Quebbemann}
}

@article{Schafer,
  author    = {R. D. Schafer},
  journal   = {Amer. J.~Math.},
  number    = {2},
  pages     = {435--446},
  title     = {On the Algebras Formed by the {C}ayley--{D}ickson Process},
  volume    = {76},
  year      = {1954},
  doi       = {10.2307/2372583}
}

@article{Viazovska8,
  author    = {Maryna Viazovska},
  title     = {The sphere packing problem in dimension $8$},
  journal   = {Annals of Mathematics},
  volume    = {185},
  number    = {3},
  pages     = {991--1015},
  year      = {2017},
  publisher = {Department of Mathematics of Princeton University},
  doi       = {10.4007/annals.2017.185.3.7}
}

@book{Voight,
  author    = {Voight, John},
  year      = {2021},
  title     = {Quaternion Algebras},
  series    = {Graduate Texts in Mathematics},
  volume    = {288},
  publisher = {Springer Nature},
  doi       = {10.1007/978-3-030-56694-4}
}

@article{Wilson,
  title     = {Octonions and the {L}eech lattice},
  author    = {Wilson, Robert A.},
  journal   = {Journal of Algebra},
  volume    = {322},
  number    = {6},
  pages     = {2186--2190},
  year      = {2009},
  publisher = {Elsevier},
  doi       = {10.1016/j.jalgebra.2009.04.015}
}

\end{document}